\documentclass[12pt]{amsart}
\usepackage{amsthm}
\usepackage{graphicx}
\usepackage{amssymb}

\newtheorem{tw}{Theorem}
\newtheorem{lem}[tw]{Lemma}
\newtheorem{wn}[tw]{Cororally}

\theoremstyle{remark}

\newcommand{\cal}[1]{\mathcal{#1}}

\newcommand{\map}[3]{#1\colon #2\to #3}
\newcommand{\Map}[2]{\map{#1}{#2}{#2}}
\newcommand{\fii}{\varphi}

\newcommand{\field}[1]{\mathbb{#1}}
\newcommand{\pp}{\field{P}}
\newcommand{\zz}{\field{Z}}
\newcommand{\rr}{\field{R}}

\newcommand{\Mod}[3]{#1\equiv #2\pmod{#3}}
\newcommand{\fal}[1]{\widetilde{#1}}
\newcommand{\gen}[1]{\langle #1 \rangle}
\newcommand{\lst}[2]{{#1}_1,\dotsc,{#1}_{#2}}

\newcommand{\SP}{\:}
\newcommand{\st}{\;|\;}

\newcommand{\Mh}{{\cal{M}}^{h}_{g}}
\newcommand{\Ms}{{\cal{M}}_{0,2g+2}}
\newcommand{\Mr}{{\cal{M}}_{0,r}}

\title[Conjugacy classes of finite\ldots]{Conjugacy classes of finite subgroups of certain mapping class groups.}

\author{Micha\l\ Stukow}

\thanks{Supported by BW5100-5-0227-2}
\thanks{{\bf Acknowledgement.} Author wishes to thank the referee for his/her helpful suggestions.}

\address[]{
Institute of Mathematics, University of Gda\'nsk, Wita Stwosza 57,
80-952 Gda\'nsk, Poland }

\email{trojkat@math.univ.gda.pl}

\keywords{Hyperelliptic mapping class group, finite subgroups of mapping class group}
\subjclass{Primary 57M60; Secondary 57N05}

\begin{document}
\begin{abstract}
We give a complete description of conjugacy classes of finite subgroups of the mapping class group
of the sphere with $r$ marked points. As a corollary we obtain a description of conjugacy classes
of maximal finite subgroups of the hyperelliptic mapping class group. In particular, we prove that
for a fixed genus $g$ there are at most five such classes.
\end{abstract}

\maketitle%

\section{Introduction}
Let $\Mr$ be the mapping class group of the sphere with $r\geq 3$ marked points, where we allow
maps to permute the set of marked points. In \cite{Buskirk} Gillette and Van Buskirk proved, using
purely algebraic methods, the following theorem.
\begin{tw}
$\Mr$ contains an element of finite order $n$ if and only if $n$
divides one of $r$, $r-1$, $r-2$.
\end{tw}
Later on the stronger version of this theorem was obtained as a by-product of certain
considerations of Harvey and Maclachlan \cite{MacHav}.
\begin{tw}
Every element of finite order in $\Mr$ is contained in a maximal
cyclic subgroup of order $r$, $r-1$ or $r-2$, and all such
subgroups of the same order are conjugate.
\end{tw}

 In this paper we extend the above results. For every finite
subgroup $N$ of $\Mr$ we find a maximal finite subgroup $M$ of $\Mr$ containing $N$. Furthermore we
give a complete description of conjugacy classes of finite subgroups of $\Mr$. As a corollary we
obtain a description of conjugacy classes of maximal finite subgroups of the hyperelliptic mapping
class group.

\section{Maximal finite subgroups of $\Mr$}
This section is devoted to the proof of the following theorem.
\begin{tw}
Finite subgroup $N$ of $\Mr$ is a maximal finite subgroup of $\Mr$ if and only if $N$ is isomorphic
to one of the following:
\begin{enumerate}
 \item the cyclic group $\zz_{r-1}$ if $r\neq 4$,
 \item the dihedral group $D_{r}$
 \item the dihedral group $D_{r-2}$ if $r=5$ or $r\geq 7$,
 \item the alternating group $A_4$ if $\Mod{r}{4\text{ or }10}{12}$
 \item the symmetric group $S_4$ if $\Mod{r}{0,2,6,8,12,14,18 \text{ or }20}{24}$
 \item the alternating group $A_5$ if $\Mod{r}{0,2,12,20,30,32,42\text{ or }50}{60}$
\end{enumerate}

\end{tw}
Let $N$  be a finite subgroup of $\Mr$. By the positive solution to the Nielsen realisation problem
\cite{Kerk1}, we could assume that $N$ is a group of automorphisms of the Riemann sphere $\pp^1$.
From the classification of such groups we have the following five possibilities:
\begin{enumerate}
 \item[1.] $N\cong\zz_n$ is generated by a rotation of order $n$,
 \item[2.] $N\cong D_n$ is generated by two rotations of order $n$ and $2$ respectively,
 \item[3.] $N\cong S_4$ is the symmetry group of a cube (octahedron),
 \item[4.] $N\cong A_5$ is the symmetry group of a dodecahedron (icosahedron), and
 \item[5.] $N\cong A_4$ is the symmetry group of a tetrahedron.
\end{enumerate}
In each case we shall find a maximal finite subgroup $M$ of $\Mr$ containing $N$. What is more,
from the below analysis it will be clear how to construct each of maximal subgroups listed in the
above theorem.

\medskip \noindent{\bf{1.}\SP}As a fundamental set for the action of a cyclic group $N$ we can
choose the interior of the bigon $F$ with vertices in the south and north poles together with one
of its edges (Figure 1).

\begin{figure}[h]
\includegraphics[width=0.8\textwidth]{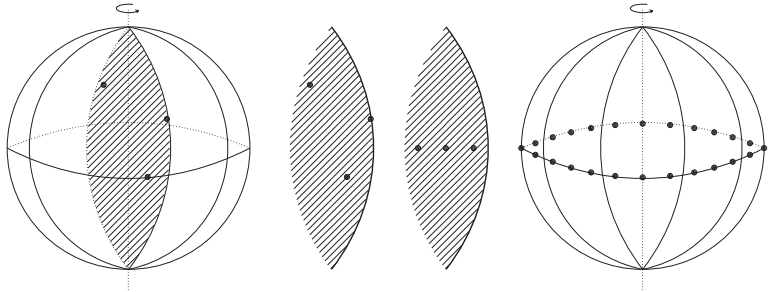}
\caption{Fundamental set for the action of $\zz_n$}
\end{figure}

Let $\cal{P}$ be the set of these marked points $P_i$, which have trivial stabiliser in $N$ and say
it has $s$ elements. Without loss of generality, we can assume that $N$ acts on the unit sphere
$S^2$ in $\rr^3$ by orthogonal rotations having the pair of poles as fixed points.

We claim, that we can also assume that $\cal{P}$ is a set of points on the equator which are
vertices of a regular $s$-gon. To prove this, we will show how to "move" points of $\cal{P}$ to
such symmetric position without changing the geometry of the action of $N$.

The canonical projection $\map{\pi}{S'}{S'/N}$, where $S'$ is obtained from $S^2$ by removing the
poles is a covering. Since $\pi(S')$ is homeomorphic to a two--punctured sphere,  we could find a
homeomorphism $\Map{f}{\pi(S')}$ isotopic to the identity which maps $\pi(\cal{P})$  onto
$\pi(\cal{P}')$, where $\cal{P}'$ is a set of $s$ points on the equator which are vertices of a
regular $s$-gon. Now $f$ has a lifting and if we extend it to $\Map{\fal{f}}{S^2}$  we see that
$\fal{f}N\fal{f}^{-1}=N$ acts on $\fal{f}(S^2)=S^2$ in such a way that $\cal{P}'$ is the set of
marked points (Figure 1). Call this model for the action of $N$ \emph{regular}.

We would like to point out that although two topological models for the action of $N$ (with the set
of marked points with trivial stabiliser equal to $\cal{P}$ and $\cal{P}'$ respectively) are models
for the action of the same subgroup $N$ of $\Mr$, they could correspond to non isomorphic
realisations of $N$ as an automorphism group of punctured Riemann sphere.

Now, using regular model for the action of $N$, it is obvious that $N$ is contained in the cyclic
group $N'$ of order~$s$.

If $s=r-2$ or $s=r$ (i.e. both poles are marked or not), $N'$ is contained in the dihedral group
$N''$ generated by $N'$ and a half turn leaving the equator invariant; we fall into the second
case. Therefore we could assume that $s=r-1$. To complete the case it is enough to prove that for
$r\neq 4$, $N'$ is a maximal subgroup of $\Mr$.

If a group $\zz_n$ acts on the sphere with $r$ marked points then at least one of them has trivial
stabiliser (since $r\geq 3$) and its orbit has length $n$, so $n\leq r$. This proves that the
action of $\zz_{r-1}$ cannot be extended to the action of any cyclic group. Similarly, if we have
the action of a dihedral group $D_n$ then the stabiliser of at least one of the marked points has
order one or two, so $n\leq r$. What is more, there is no action by the group $D_{r-1}$.

The only cyclic subgroups of $A_4,S_4,A_5$ are, respectively, $\zz_2,\zz_3$ and $\zz_2,\zz_3,\zz_4$
and $\zz_2,\zz_3,\zz_5$, so we could restrict ourselves to the case $r\leq 6$. The maximal order of
a stabiliser in $A_4$ is $3$, so the minimal length of an orbit is $4$ and we see that the action
of $A_4$ is possible only for $r\geq 4$. Similarly, actions of $S_4$ and $A_5$ are possible only
for $r\geq 6$ and $r\geq 12$ respectively. Therefore the possible extension of the action of
$\zz_{r-1}$ could exits only for $r=4$. Conversely if $r=4$, using the above technique of obtaining
regular action, we could extend action of $\zz_3$ to the action of $A_4$.

\medskip
\noindent {\bf{2.\SP}} If $N$ is a dihedral group, then using similar methods such as in the case
of cyclic groups we could construct a regular model for the action of $N$, i.e. such that the
action of the canonical cyclic subgroup $K$ of $N$ is regular and $N$ is generated by $K$ and a
half turn, leaving the equator invariant.

Now it is clear that, depending on whether the poles are marked or not, this group is contained in
${\rm D}_{r-2}$ or ${\rm D}_r$. Using similar arguments like in the case of a cyclic group, we
prove that $D_r$ for $r\geq 3$ and $D_{r-2}$ for $r=5$ or $r\geq 7$ are maximal. If $r=3$ then
$D_{r-2}$ is a cyclic group and if $r=4$ or $r=6$ then we could extend the action of $D_{r-2}$ to
the action of $D_4$ and $S_4$ respectively.

\medskip
\noindent{\bf{3.\SP}}If $N$ is the group of orientation preserving symmetries of a cube, then we
have three orbits of points with nontrivial stabiliser in $N$: the centres of faces, the centres of
edges and the vertices with stabilisers of order $4,2$ and $3$, respectively. Since the length of
an orbit is equal to the index of the stabiliser, these orbits have lengths $6,12$ and $8$,
respectively.

Let $k$ be  the number of marked points in a fixed fundamental set $F$, which have trivial
stabiliser. Suppose that centres of faces, edges and vertices are marked points. Then counting
lengths of orbits, we have $r=24k+8+12+6$, i.e. $k=\frac{r-26}{24}$. In similar way, according to
whether centres of faces, edges or vertices are marked points or not, we obtain other possible
values of $k$ (Table 1).

\begin{table}
\begin{tabular}{|c||c|c|c||c|c|c|}
\hline
&F&E&V&cube&dodecahedron&tetrahedron\\
\hline
(a)&--&--&--&$r/24$&$r/60$&$r/12$\\

(b)&+&--&--&$(r-6)/24$&$(r-12)/60$&$(r-4)/12$\\

(c)&--&+&--&$(r-12)/24$&$(r-30)/60$&$(r-6)/12$\\

(d)&+&+&--&$(r-18)/24$&$(r-42)/60$&$(r-10)/12$\\

(e)&--&--&+&$(r-8)/24$&$(r-20)/60$&$(r-4)/12$\\

(f)&+&--&+&$(r-14)/24$&$(r-32)/60$&$(r-8)/12$\\

(g)&--&+&+&$(r-20)/24$&$(r-50)/60$&$(r-10)/12$\\

(h)&+&+&+&$(r-26)/24$&$(r-62)/60$&$(r-14)/12$\\
\hline
\end{tabular}
\medskip
\caption{Possible values of $k$}
\end{table}

%To obtain that table we proceed as follows. Suppose, for example that centres of faces, edges or
%vertices are marked points (case (h)). Counting lengths of orbits, we have $r=24k+8+12+6$, i.e.
%$k=\frac{r-26}{24}$.

Observe that for the fixed number $r$ of distinguished points, only one of the above cases could
occur. Since there is no monomorphism of $S_4$ into $\zz_n,D_n,A_4,A_5$, every such group is maximal.%

\medskip
\noindent {\bf4.\SP}The case of the group of orientation preserving symmetries of a dodecahedron is
fully analogous to the case of a cube and we obtain also eight possibilities (Table~1).

\medskip
\noindent
 {\bf 5.\SP}If $N$ is the group of rotations of a
regular tetrahedron, then by the same reasons as before, we have eight possibilities (Table~1).

Observe that, since centres of faces are vertices of the dual tetrahedron, it is not difficult to
show that groups in cases (b), (e) and (d), (g) are conjugate.

However the problem of maximality here is more involved due to the fact that there are natural
monomorphisms of $A_4$ into $S_4$ and $A_5$ induced by embeddings of a tetrahedron into a cube and
dodecahedron (Figure 2).

\begin{figure}[h]%\begin{center}
\includegraphics[width=0.6\textwidth]{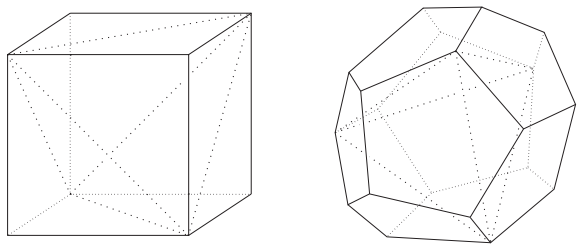} %[width=1\textwidth]
\caption{Embeddings of a tetrahedron}
\end{figure}
Suppose that $N$ has an extension $N'\cong S_4$. Since $S_4$ contains only one subgroup isomorphic
to $A_4$, if we realise $N'$ as the group of rotations of a cube, we see that such extension could
exist only when vertices and centres of faces of the tetrahedron (defined as points of $S^2$ with
stabiliser in $N$ of order $3$) are of the same type (marked or not), because after extension they
fall into one orbit. In particular in cases (b), (d), (e), (g) such extension could not exist.

The same conclusion is true for the extension of $N$ to the action of $A_5$, however now one has to
be more careful, because $A_5$ contains five subgroups isomorphic to $A_4$.

Conversely, using the same method as in the case of a cyclic group, it is easy to obtain in cases
(a), (c), (f), (h) extension of $N$ to the action of $S_4$.

\section{Conjugacy classes of finite subgroups of $\Mr$}
Let $N$ be a finite subgroup of $\Mr$. From the previous section, follows that $N$ is isomorphic
either to the cyclic group $\zz_n$ or to the dihedral group $D_n$ or to the symmetry group of a
cube (octahedron) $S_4$ or the symmetry group of a dodecahedron (icosahedron) $A_5$ or else  the
symmetry group of a tetrahedron $A_4$. We are now ready to extend this description to conjugacy
classes.
\begin{tw}
The set of conjugacy classes of subgroups isomorphic to $N$ has exactly  two elements if $N\cong
\zz_2$ with  $r$ even or  $N\cong D_n$ with  $2n|r$ or $2n|r-2$ and one element for the remaining
cases.
\end{tw}
A necessary condition for two subgroups $N_1,N_2$  of $\Mr$ to be conjugate is that corresponding
orbits of points with nontrivial stabiliser are of the same type, i.e. they are marked or not.
Counting for each possible type of the action the number of points of  a fundamental set with 
trivial stabiliser one can see  that for most cases different actions live on spheres with
different numbers of distinguished points. The only exceptions are those of the first part of the
theorem.

If $r$ is even, then according to whether poles are marked or not, we have two non-conjugate
actions of $\zz_2$. Now suppose $2n|r$. If we have an action of $D_n$ then there are two
possibilities to distribute $r$ marked points: the first is that all marked points form one orbit
of length $2n$, and the second when marked points form two orbits of order $n$. This two
possibilities give two non-conjugate actions of $D_n$. Similar situation holds for $2n|r-2$.

Conversely, let $N_1$ and $N_2$  be two isomorphic subgroups of $\Mr$ with the same structure of
orbits of points with nontrivial stabiliser. From the previous section it follows that there exist
homeomorphisms $\map{f_i}{X}{X_i}$ for $i=1,2$, such that $N_i'=f_iN_if_i^{-1}$ acts on the sphere
$X_i$ in a "regular" way (i.e. $X_1=X_2=S^2$, sets of marked points coincide, both groups act by
orthogonal rotations and as a groups of isometries of $S^2$, $N_1=N_2$). Let $\map{\fii}{X_1}{X_2}$
be the identity map. Observe that $\fii$ is the identity as a map of $S^2$, though as an element of
$\Mr$ it could permute marked points. So for $\psi=f_1^{-1}\fii^{-1}f_2$ we have
\[\psi N_2 \psi^{-1}=f_1^{-1}\fii^{-1}f_2 N_2 f_2^{-1}\fii f_1=f_1^{-1}\fii^{-1}N_2'\fii f_1
=f_1^{-1}N_1'f_1=N_1\]%
\begin{wn}
Two maximal finite subgroups of $\Mr$ are conjugate if and only if they are isomorphic.
\end{wn}
\section{Conjugacy classes of maximal finite subgroups of $\Mh$}
Suppose that a closed, connected orientable surface $T_g$ of genus $g\geq 2$ is embedded in $\rr^3$
in such a manner that it is invariant under the half turn $\rho$ about the $y$-axis (figure 3).

\begin{figure}[h]
\includegraphics[width=0.7\textwidth]{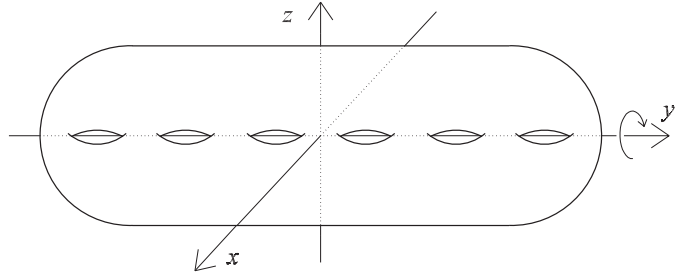} %[width=1\textwidth]
\caption{Hyperelliptic involution $\rho$}
\end{figure}
The hyperelliptic mapping class group $\Mh$ is defined to be the centraliser of $\rho$ in the
mapping class group of $T_g$. By \cite{Bir-Hil} the quotient $\Mh/\gen{\rho}$ is isomorphic to the
mapping class group $\Ms$ of a sphere with $2g+2$ marked points $\lst{P}{2g+2}$. For simplicity we
will identify $\Mh/\gen{\rho}$ with $\Ms$. Let $\map{\pi}{\Mh}{\Ms}$ be the canonical projection
and ${\cal{K}},{\cal{L}}$ be sets of maximal finite subgroups of $\Mh$ and $\Mr$ respectively.
\begin{lem}
The map $\map{\fii}{{\cal{K}}}{{\cal{L}}}$ given by $\fii(K)=\pi(K)=K/\gen{\rho}$ is a bijection.
\end{lem}
\begin{proof}
Observe that since $\rho$ is central in  $\Mh$, it is contained in each maximal finite subgroup of
$\Mh$, so $\fii$ is well defined. Clearly $\fii$ is surjective, so let us prove that it is also
injective. If $\pi(K)=L$ for $K\in\cal{K}$ and $L\in\cal{L}$ then $K\leq\pi^{-1}(L)$. Comparing
orders, we obtain $K=\pi^{-1}(L)$.
\end{proof}
\begin{lem}
Under above assumptions, $\fii$ induces a bijection between conjugacy classes of subgroups in
${\cal{K}}$ and ${\cal{L}}$.
\end{lem}
\begin{proof}
Obviously conjugate elements of $\cal{K}$ are mapped onto conjugated elements of $\cal{L}$. The
converse statement, follows from bijectivity of~$\fii$.
\end{proof}

Therefore, using above lemmas and known description of liftings of finite subgroups of $\Ms$ to
$\Mh$ \cite{HyperGG}, we could reformulate theorems of previous sections in terms of the group
$\Mh$.
\begin{tw}
Finite subgroup $N$ of $\Mh$ is a maximal finite subgroup of $\Mh$ if and only if $N$ is isomorphic
to one of the following:
\begin{enumerate}
\item $\zz_{4g+2}$
\item $V_{2g+2}=\gen{x,y\st x^4,y^{2g+2},(xy)^2,(x^{-1}y)^2}$
\item $U_{2g}=\gen{x,y\st x^2,y^{4g},xyxy^{2g+1}}$ if $g\geq 3$
\item $\zz_2\times A_4$ if $\Mod{g}{1}{6}$
\item $SL(2,3)=
\gen{x,y\st x^4,y^3,(xy)^3,yx^2y^{-1}x^2}$ if $\Mod{g}{4}{6}$
\item $\zz_2\times S_4$ if $\Mod{g}{3\text{ or }11}{12}$
\item $W_1=\gen{x,y\st x^2,y^3,(xy)^4(yx)^4,(xy)^8}$ if $\Mod{g}{2
\text{ or }6}{12}$
\item $W_2=\gen{x,y\st x^4,y^3,yx^2y^{-1}x^2,(xy)^4}$ if $\Mod{g}{5
\text{ or }9}{12}$
\item $W_3=\gen{x,y\st x^4,y^3,(xy)^8,x^2(xy)^4}$ if $\Mod{g}{0
\text{ or }8}{12}$
\item $\zz_2\times A_5$ if $\Mod{g}{5,9,15\text{ or }29}{30}$
\item $SL(2,5)=\gen{x,y\st x^4,y^3,(xy)^5,yx^2y^{-1}x^2}$ if $\Mod{g}{0,14,20\text{ or }24}{30}$
\end{enumerate}
Two maximal finite subgroups of $\Mh$ are conjugate within $\Mh$ if and only if they are
isomorphic.
\end{tw}
Groups in the first part of the theorem are respectively liftings of (1) $\zz_{r-1}$, (2) $D_r$,
(3) $D_{r-2}$, (4)--(5) $A_4$, (6)--(9) $S_4$, (10)--(11) $A_5$, where $r=2g+2$. The reason for the
assumption $g\geq 3$ in case (3) is that for $g=2$ the underlying group $D_{r-2}$ is not maximal.
\begin{wn}
If $\Mod{g}{0,5,9,14,15,20,24 \text{ or } 29}{30}$ then $\Mh$ contains exactly five conjugacy
classes of maximal finite subgroups. For $g=2$ there are three and for other values of $g$ there
are exactly four such classes.
\end{wn}

\end{document}